%% file: asymptotic.shade.ub.tex
\documentclass{amsart}

\usepackage{enumitem}
\usepackage{hyperref}
\usepackage{calc}

\makeatletter
\include{macros}

\newtheorem{thm}{Theorem}

\newtheorem{lem}{Lemma}
\newtheorem{cor}[lem]{Corollary}
\newtheorem{prop}[lem]{Proposition}

\newtheorem{conjecture}{Conjecture}
\newtheorem{emp}{\empstring}

\theoremstyle{definition}

\theoremstyle{remark}

\newcommand{\empstring}{Dichotomy}

\setenumerate{leftmargin=*, label=\tu{(\alph*)}, ref=\alph*, widest=b}

\begin{document}

\title[Asymptotic upper bounds on shades of $t$-intersecting
  families]{Asymptotic upper bounds on the shades\\ of $t$-intersecting
  families}
\author{James Hirschorn}
\address{Thornhill, ON, Canada}
\email{\href{mailto:j_hirschorn@yahoo.com}{j\textunderscore hirschorn@yahoo.com}}
\urladdr{\href{http://www.logic.univie.ac.at/~hirschor}{http://www.logic.univie.ac.at/\textasciitilde  hirschor}}
\date{August 10, 2008}
\maketitle

\begin{abstract}
We examine the $m$-shades of $t$-intersecting families of $k$-subsets
of $[n]$, and conjecture on the optimal upper bound on their
cardinalities. This conjecture extends Frankl's General Conjecture
that was proven true by Ahlswede--Khachatrian. From this we deduce the
precise asymptotic upper bounds on the cardinalities of $m$-shades of
$t(m)$-intersecting families of $k(m)$-subsets of $[2m]$, as
$m\to\infty$. A generalization to cross-$t$-intersecting families is
also considered.
\end{abstract}

\section{Introduction}
\label{sec:motivation}

The paper~\cite{Hirs1} was concerned with the dichotomy below of
descriptive set theory. This dichotomy is aimed towards research on
a fundamental question of set theoretic forcing, of whether Cohen and
random forcing together form a basis for all nontrivial Souslin
(i.e.~``simply'' definable) ccc (i.e.~no uncountable antichains) 
posets (asked by Shelah in~\cite{MR1303493}).

\begin{emp}
\label{u-1}
Every analytic \tu(i.e.~projection of a closed subset of the
\tu{``}plane\tu{'')} family $\A$ of infinitely branching subtrees of $\twoseq$
satisfies at least one of the following\textup:
\begin{enumerate}
\item\label{item:2} There exists a colouring $c:\twoseq\to\two$ and a
  such that $S(n)$ is nonhomogeneous for $c$ for all but finitely many
  $n\in\N$, for every $S\in\A$.
\item\label{item:1} The poset $(\A,\subseteq)$ has an uncountable antichain.
\end{enumerate}
\end{emp}

It turned out that obtaining tight upper bounds on the $m$-shades of $t$-intersecting
families of $k$-subsets of $[m]$, was relevant to
dichotomy~\ref{u-1}. The connection is described in lemma~\ref{l-2}
below (cf.~\cite{Hirs1} for details).

\subsection{Shades}
\label{sec:shades}

One of the basic notions in Sperner theory is the \emph{shade} (also
called \emph{upper shadow}) of a set or a family
of sets (see e.g.~\cite{MR2003b:05001},\cite{MR1429390}). 
For a subset $x$ of a fixed set~$S$, the shade of $x$ is
\begin{equation}
  \label{eq:8}
  \nabla(x)=\{y\subseteq S:x\subset y\text{ and }|y|=|x|+1\},
\end{equation}
and the shade of a family~$X$ of subsets of $S$ is 
\begin{equation}
\label{eq:11}
\nabla(X)=\bigcup_{x\in X}\nabla(x).
\end{equation}
Recall that the \emph{$m$-shade} (also called \emph{upper $m$-shadow}
or \emph{shade at the $m\Th$ level}) of~$x$ is
\begin{equation}
  \label{eq:13}
  \shadeto m(x)=\{y\subseteq S:x\subseteq y\text{ and }|y|=m\},
\end{equation}
and $\shadeto m(X)=\bigcup_{x\in X}\shadeto m(x)$. We follow the
Sperner theoretic conventions of writing $[m,n]$ for the set
$\{m,m+1,\dots,n\}$ and  $[n]$ for the set $[1,n]=\{1,\dots,n\}$. 

We introduce the following notation for colouring sets with two colours. For a set
$S$, let $S \choose [m]$ denote the collection of all colourings $c:S\to\two$ with
$|c\inv(0)|=m$, i.e.~$c\inv(0)=\{j\in S:c(j)=0\}$. This is related to shades,
because for all $c\in{S\choose[m]}$ and all $x\subseteq S$,
\begin{multline}
  \label{eq:49}
  x\text{ is homogeneous for }c \Iff c\inv(0)\in\shadeto m(x)\Or c\inv(1)\in\shadeto{|S|-m}(x).
\end{multline}

When a nonhomogeneous colouring is desired, it is most efficient to use colorings in
$S\choose [m]$ for $|S|=2m$. Equation~\eqref{eq:49} immediately gives us:

\begin{lem}
\label{l-2}
Suppose $X$ is a family of subsets of $[2m]$. Then
\begin{equation*}
  \left|\left\{c\in{[2m]\choose[m]}:\exists x\in X
      \spc x\textup{ is homogeneous for }c\right\}\right|
  \le2|\shadeto m(X)|
\end{equation*}
\textup(the shades are with respect to $S=[2m]$\textup).
\end{lem}

\section{Upper bounds}
\label{sec:upper-bounds}

Recall that a family $\A$ of sets is \emph{$t$-intersecting} if $|E\cap F|\ge t$ for all
$E,F\in\A$; and a pair~$(\A,\B)$ of families of subsets of some fixed 
set are \emph{cross-$t$-intersecting} if 
\begin{equation}
  \label{eq:26}
  |E\cap F|\ge t\espc\text{for all $E\in\A$, $F\in\B$}.
\end{equation}
Thus $\A$ is $t$-intersecting iff $(\A,\A)$ is cross-$t$-intersecting.

We use the standard notation $S\choose k$ to denote the collection of
all $k$-subsets of~$S$, and hence $[n]\choose k$ denotes the collection
of all subsets of~$[n]$ of cardinality $k$. Let $I(n,k,t)$ denote
the family of all $t$-intersecting subfamilies of $[n]\choose k$
(where $t\le k\le n$). 
Define the function
\begin{equation}
  \label{eq:31}
  M(n,k,t)=\max_{\A\in I(n,k,t)}|\A|.
\end{equation}
The investigation into the function $M$ and the structure of the
maximal families was initiated by Erd\H os--Ko--Rado in 1938, but not
published until~\cite{MR0140419}. In this paper, they gave a complete
solution for the case $t=1$, and posed what became one of the most famous
open problems in this area. The following so called
\emph{$4m$-conjecture} for the case $t=2$:
\begin{equation}
  \label{eq:39}
  M(4m,2m,2)=\frac 12 \biggl({{4m}\choose{2m}}-{{2m}\choose m}^{\hspb 2}\biggr).
\end{equation}
We briefly explain the significance of the right hand side expression.
Define families
\begin{equation}
  \label{eq:62}
  \F_i(n,k,t)=\left\{F\in{[n]\choose k}:|F\cap[t+2i]|\ge t+i\right\}
  \espc\text{for $0\le i\le\frac{n-t}2$}.
\end{equation}
Clearly each $\F_i(n,k,t)$ is $t$-intersecting. In the special case
where $n=2k=2m$ and $t=2s$,
we can easily compute the cardinality of the corresponding
$\F_i$ using the fact that $[2m]\setminus F$ is an $m$-set for all $F\in{[2m]\choose m}$,
i.e.~$|\F_i(2m,m,2s)|$ equals
\begin{equation}
  \label{eq:68}
  \frac12\left({2m\choose m}-\sum_{j=-(s-1)}^{s-1}
    {2(s+i)\choose s+i+j}{2m-2(s+i)\choose m-(s+i+j)}\right).
\end{equation}
Then plugging in $m:=2m$ (i.e.~$2m$ for $m$), $s:=1$ and $i:=m-1$ we see
that the right hand side of 
equation~\eqref{eq:39} is equal to the cardinality of
$\F_{m-1}(4m,2m,2)$. The $4m$-conjecture was
generalized by Frankl in 1978 (\cite{MR519277}) as follows: For all
$1\le t\le k\le n$,
\begin{equation}
  \label{eq:63}
  M(n,k,t)=\max_{0\le i\le\frac{n-t}2}|\F_i(n,k,t)|.
\end{equation}
In 1995, the general conjecture was proven true by
Ahlswede--Khachatrian in\linebreak \cite{MR1429238}, where they moreover
established that the optimal families in $I(n,k,t)$ are equal to one
of the families $\F_i(n,k,t)$ up to a permutation of $[n]$. This finally
settled the $4m$-conjecture, and moreover proved that the maximal family in
$I(4m,2m,2)$ is isomorphic to $\F_{m-1}(4m,2m,2)$. 

For reasons alluded to in lemma~\ref{l-2}, 
it is upper bounds on the cardinality of the shades of $t$-intersecting families that we
are interested in, rather than upper bounds on the families
themselves. While there are numerous results giving lower bounds on
the size of shadows/shades, upper bounds seem to be rather scarce.
Perhaps this is because they are not very good. For example, the
following is from~\cite{MR1044236}, where $2^{S}$ denotes the power
set of $S$. 

\begin{thm}[Kostochka, 1989]
Suppose that $\A\subseteq 2^{[n]}$ is a Sperner family.
Then $\nabla(\A)\le0.724\cdot2^n$. 
\end{thm}

\noindent Moreover, the best upper bound is
known to be greater than $0.5\cdot2^n$. 

However, in the case of $t$-intersecting families, the shade is also
$t$-intersecting. Define for $1\le t\le k\le m\le n$,
\begin{equation}
  \label{eq:64}
  M_0(n,m,k,t)=\max_{\A\in I(n,k,t)}|\shadeto m(\A)|,
\end{equation}
i.e.~$M_0(n,m,k,t)$ is the maximum size of the $m$-shade of a
$t$-intersecting family of $k$-subsets of $[n]$. Thus we have
\begin{equation}
  \label{eq:65}
  M_0(n,m,k,t)\le M(n,m,t),
\end{equation}
but this is not optimal. Indeed we make the following easy observations.

\begin{lem}
\label{l-7}
For all $1\le t\le k\le m\le n$,
\begin{enumerate}
\item\label{item:12} $\F_i(n,k,t)=\emptyset$ for all $i>k-t$,
\item\label{item:13} $\shadeto m(\F_i(n,k,t))=\F_i(n,m,t)$ for all
  $0\le i \le \min\bigl(k-t,\frac{n-t}2\bigr)$.  
\end{enumerate}
\end{lem}

This leads to the following conjecture.

\begin{conjecture}
\label{j-1}
$\displaystyle M_0(n,m,k,t)=\max_{0\le i\le\min(k-t,\frac{n-t}2)}|\F_i(n,m,t)|$.
\end{conjecture}

\noindent Note that conjecture~\ref{j-1} is correct so long as the
optimal families that we are taking the $m$-shades of are among the
$\F_i(n,k,t)$.

\subsection{Asymptotic behaviour}
\label{sec:asymptotic-behaviour}

Not surprisingly, for the purpose of our set theoretic dichotomy we were
interested in the asymptotic behaviour of the upper bounds, i.e.~as
$n\to\infty$. Furthermore, we were interested in the $m$-shade of
subsets of $[2m]$, i.e.~$n=2m$. We would like something to
the effect that the maximum proportion of the $m$-shade to the entire family
$[2m]\choose m$ goes to $0$; symbolically, 
\begin{equation}
\lim_{\substack{m\to\infty\\ t(m)\le k(m)\le m}}
\frac{M_0(2m,m,k(m),t(m))}{{{2m} \choose m}}=0.\label{eq:66}
\end{equation}
However, this is false, because for example the optimal family for the
$4m$-conjecture (cf.~equation~\eqref{eq:39}) gives us
$\lim_{m\to\infty}M_0(2m,m,m,2)\div {2m\choose m}=\frac12$.
In fact, equation~\eqref{eq:66} fails whenever $t(m)$ is bounded:
\begin{equation}
  \label{eq:69}
  \lim_{m\to\infty}\frac{M_0(2m,m,t,t)}{{2m\choose m}}\ge\frac1{2^t}.
\end{equation}
This can be seen by noting that $\shadeto
m(\F_0(2m,t,t))=\F_0(2m,m,t)$, 
and that $|\F_0(2m,m,t)|={2m-t\choose
  m-t}$. 

Moreover, as we shall demonstrate, even if $t(m)\xrightarrow{m\to\infty}\infty$,
equation~\eqref{eq:66} can still fail if $t(m)$ is too small compared
with $k(m)$. Henceforth, we shall make the simplification that
$k(m)=o(m)$, i.e.~$\lim_{m\to\infty}k(m)\div m=0$; this was the only case used in our
application. 

Before proceeding further, 
recall the de Moivre--Laplace theorem (cf.~\cite{63.1069.01}), roughly stating that
the binomial series of $(p+q)^n$ has most of the sum
concentrated in the order of~$\sqrt n$ terms around the center: For
all $0\le a,b<\infty$,
\begin{equation}
  \label{eq:77}
  \lim_{n\to\infty}\sum_{j=- \lfloor a\sqrt{n\div 2}\rfloor}^{\lfloor
    b\sqrt{n\div 2}\rfloor}
  \frac{{2n\choose n+j}}{4^n}
  =\Phi(b)-\Phi(-a),
\end{equation}
where
\begin{equation}
  \label{eq:84}
  \Phi(t)=\frac1{\sqrt{2\pi}}\int_{-\infty}^{t}e^{-{x^2}\div 2}\der x
\end{equation}
is the cumulative distribution function of the standard normal
distribution. Recall that
\begin{equation}
  \label{eq:85}
  \Phi(t)-\Phi(-t)=2\Phi(t)-1.
\end{equation}
We write $f\sim g$ to indicate asymptotic equality,
i.e.~$\lim_{n\to\infty}\frac{f(n)}{g(n)}=1$. 
The theorem also tells us that
\begin{equation}
  \label{eq:81}
  \frac{{2n\choose n+j}}{4^n}
  \sim\frac{e^{-(j\div\sqrt{n\div2})^2\div2}}{\sqrt{\pi n}},
\end{equation}
and moreover that the convergence is uniform over $j\div\sqrt{n\div2}$
in the range $[a,b]$. 

We derive a version of~\eqref{eq:77} for the identity 
${2n\choose n}=\sum_{j=-k}^k{2k\choose k+j}{2(n-k)\choose (n-k)-j}$.

\begin{lem}
\label{l-3}
Assume $k(n)=o(n)$ and $\lim_{n\to\infty}k(n)=\infty$. Then
\begin{equation}
\label{eq:20}
\lim_{n\to\infty}\sum_{j=-\lfloor a\sqrt{k(n)\div 2}\rfloor}
^{\lfloor b\sqrt{k(n)\div 2}\rfloor}\frac{{2k(n)\choose k+j}
  {2(n-k(n))\choose (n-k(n))-j}}{{2n\choose n}}=\Phi(b)-\Phi(a).
\end{equation}
\end{lem}
\begin{proof}
By equation~\eqref{eq:81} and the assumptions on $k$, 
we have ${2k(n)\choose k(n)+j}\sim\frac{4^{k(n)}}{\sqrt{\pi k(n)}}\cdot
e^{-(j\div\sqrt{k(n)\div 2})^2\div2}$ and 
${2(n-k(n))\choose(n-k(n))-j}\sim
\frac{4^{n-k(n)}}{\sqrt{\pi(n-k(n))}}\cdot
e^{-(j\div\sqrt{(n-k(n))\div2})^2\div 2}$, with uniform convergence for 
$j\div{\sqrt{k(n)\div2}}\in[a,b]$. Therefore,
\begin{equation}
  \label{eq:28}
  {2k(n)\choose k(n)+j}{2(n-k(n))\choose(n-k(n))-j}\sim
  \frac{4^n e^{-(j\div{\sqrt{k(n)\div2}})^2\div2}}{\pi\sqrt{k(n)(n-k(n))}}
\end{equation}
with uniform convergence. Changing variables then gives
\begin{equation}
  \label{eq:29}
  \sum_{j=-\lfloor a\sqrt{k(n)\div 2}\rfloor}
^{\lfloor b\sqrt{k(n)\div 2}\rfloor}{2k(n)\choose k(n)+j}{2(n-k(n))\choose(n-k(n))-j}\sim
\frac{4^n}{\pi \sqrt{2n}}\int_a^b e^{-x^2\div2}\der x,
\end{equation}
and the right hand side expression is equal to $\frac{4^n}{\sqrt{\pi
    n}}(\Phi(b)-\Phi(a))$ which is asymptotically equal to ${2n\choose
  n}(\Phi(b)-\Phi(a))$ by equation~\eqref{eq:81} with $j=0$, as required. 
\end{proof}

\begin{lem}
\label{l-10}
Let $c>0$. Supposing $k(m)=o(m)$ and $\lim_{m\to\infty}k(m)=\infty$,
\begin{equation}
  \label{eq:78}
  \lim_{m\to\infty}\frac{|\F_{k(m)}(2m,m,c\sqrt{k(m)})|}{{2m\choose m}}
  =\frac1{\sqrt{2\pi}}\int_{c\div\sqrt2}^\infty e^{-x^2\div2}\der x.
\end{equation}
\end{lem}
\begin{proof}
Setting $s(m):=c\sqrt{k(m)}\div2$ and $i(m):=k(m)$, since
$\sqrt{s(m)+i(m)}=\sqrt{k(m)}+o(\sqrt{k(m)})$,
by lemma~\ref{l-3} and~\eqref{eq:85}
\begin{equation}
  \label{eq:83}
  \lim_{m\to\infty}\frac{\sum_{j=-(s(m)-1)}^{s(m)-1}
    {2(s(m)+i(m))\choose s(m)+i(m)+j}{2m-2(s(m)+i(m))\choose
      m-(s(m)+i(m)+j)}}{{2m\choose m}}=2\Phi(c\div\sqrt2)-1.
\end{equation}
Hence by equation~\eqref{eq:68}, the limit in~\eqref{eq:78} is equal
to $\frac12\bigl(1-(2\Phi(c\div\sqrt2)-1)\bigr)=1-\Phi(c\div\sqrt2)$
as required.
\end{proof}

For an infinite $A\subseteq\N$ we let $e_A:\N\to\N$ denote the strictly
increasing enumeration of $A$. 

\begin{cor}
\label{l-9}
Assume $k(m)=o(m)$.
Suppose there exists $c>0$ such that
\begin{equation}
  \label{eq:76}
  t(m)\le c\sqrt{k(m)}\espc\tu{for infinitely many $m$}.
\end{equation}
Then $\limsup_{m\to\infty}M_0(2m,m,k(m),t(m))\div{{2m\choose m}}>0$.
\end{cor}
\begin{proof}
Let $A=\{m\in\N:t(m)\le c\sqrt{k(m)}\}$, which is infinite by assumption. 
If $\lim_{m\to\infty}k\circ e_A(m)\ne\infty$, 
then $t(m)$ is bounded on an infinite subset $B\subseteq A$,
and the argument in equation~\eqref{eq:69} shows that
$M_0(2m,m,k(m),t(m))\div{2m\choose m}$ is bounded away from $0$ on
$B$. Otherwise, we can choose $k'$ so that $k'(m)=k(m)$ for all $m\in
A$, $\lim_{m\to\infty}k'(m)=\infty$ and $k'(m)=o(m)$.

Put $k''(m)=\max(0,k'(m)-c\sqrt{k'(m)})$.
Obviously $\lim_{m\to\infty}k''(m)=\infty$ and $k''(m)=o(m)$.
And for all $m\in A$ with $k''(m)>0$,
$\F_{k''(m)}(2m,k(m),c\sqrt{k'(m)})$ is a $t(m)$-intersecting family of
$k(m)$-subsets of $[2m]$, and furthermore its $m$-shade is
$\F_{k''(m)}(2m,m,c\sqrt{k'(m)})$ by lemma~\ref{l-7}. 
The result thus follows from lemma~\ref{l-10} with $k:=k''$.
\end{proof}

Avoiding the example of corollary~\ref{l-9},
we arrive at the following optimal conjecture, i.e.~there is no room
for improvement on $t(m)$ by corollary~\ref{l-9}.

\begin{conjecture}
\label{j-2}
$\displaystyle{\lim_{m\to\infty}\frac{M_0(2m,m,k(m),t(m))}
  {{2m\choose m}}=0}$,\newline 
whenever $k(m)=o(m)$, $\lim_{m\to\infty}k(m)=\infty$ and
$\lim_{m\to\infty}\frac{t(m)}{\sqrt{k(m)}}=\infty$.
\end{conjecture}

Let us show that conjecture~\ref{j-2} is a consequence of conjecture~\ref{j-1}.

\begin{lem}
\label{l-12}
Assume $t(m),l(m)=o(m)$ and $\lim_{m\to\infty}t(m)+l(m)=\infty$. 
Supposing that $\lim_{m\to\infty}\frac{t(m)}{\sqrt{l(m)}}=\infty$, 
\begin{equation}
  \label{eq:87}
  \lim_{m\to\infty}\frac{|\F_{l(m)}(2m,m,t(m))|}{{2m\choose m}}=0.
\end{equation}
\end{lem}
\begin{proof}
Setting $s(m):=t(m)\div2$ and $i(m):=l(m)$, 
by our assumptions on $t$ and~$l$, lemma~\ref{l-3} applies yielding
\begin{equation}
  \label{eq:88}
  \lim_{m\to\infty}\frac{\sum_{j=-(s(m)-1)}^{s(m)-1}
    {2(s(m)+i(m))\choose s(m)+i(m)+j}{2m-2(s(m)+i(m))\choose
      m-(s(m)+i(m)+j)}}{{2m\choose m}}=1.
\end{equation}
The result now follows from equation~\eqref{eq:68}. 
\end{proof}

\begin{cor}
\label{l-11}
If conjecture~\tu{\ref{j-1}} is correct then so is conjecture~\tu{\ref{j-2}}.
\end{cor}
\begin{proof}
Assume $k(m)=o(m)$, $k(m)\to\infty$ and
$\lim_{m\to\infty}t(m)\div\sqrt{k(m)}=\infty$; note that we are
implicitly assuming that $t(m)\le k(m)$ for all $m$ (so that
$M_0(2m,\allowbreak m,\allowbreak k(m),\allowbreak t(m))$ makes sense). 
Define $k'$ so that for all $m$, $0\le k'(m)\le k(m)-t(m)$ and
\begin{equation}
  \label{eq:94}
  \max_{0\le i\le k(m)-t(m)}|\F_i(2m,m,t(m))|=|\F_{k'(m)}(2m,m,t(m))|.
\end{equation}
Setting $l(m):=k'(m)$, it is clear that $t$ and $l$ satisfy the
hypotheses of lemma~\ref{l-12}. Now
\begin{equation}
  \label{eq:97}
  \lim_{m\to\infty}\frac{M_0(2m,m,k(m),t(m))}{{2m\choose m}}
  =\lim_{m\to\infty}\frac{|\F_{l(m)}(2m,m,t(m))|}{{2m\choose m}}=0,
\end{equation}
where the first equality is by conjecture~\ref{j-1}, and 
the second by lemma~\ref{l-12}. 
\end{proof}

\section{Cross-$t$-intersecting families}
\label{sec:cross-t-intersecting}

It turns out that for our application we needed upper bounds on the size
of shades of cross-$t$-intersecting families (cf.~equation~\eqref{eq:26}).
Let $C(n,k,l,t)$ be the collection of all pairs $(\A,\B)$ of cross-$t$-intersecting
families, where $\A\subseteq{[n]\choose k}$ and
$\B\subseteq{[n]\choose l}$. Then the cross-$t$-intersecting function
corresponding to $M$ is defined by
\begin{equation}
  \label{eq:57}
  N(n,k,l,t)=\max_{(\A,\B)\in C(n,k,l,t)}|\A|\cdot|\B|.
\end{equation}

There are a number of results on cross-$t$-intersecting families in
the literature; however, the state of knowledge seems very meager
compared with $t$-intersecting families. 
The following theorem, proved in~\cite{MR90g:05008}, is the strongest
result of its kind that we were able to find. 

\begin{thm}[Matsumoto--Tokushige, 1989]
$N(n,k,l,1)={n-1\choose k-1}{n-1\choose l-1}$ whenever $2k,2l\le n$.
\end{thm}

Note that this corresponds to case $t=1$ of the Erd\H os--Ko--Rado
Theorem, proved back in 1938. It is also conjectured that the EKR
Theorem does generalize:

\begin{conjecture}
\label{j-3}
$N(n,k,l,t)={n-t\choose k-t}{n-t\choose l-t}$ for all $n\ge n_0(k,l,t)$.
\end{conjecture}

Generalizing the families $\F_i$, we define
\begin{equation}
  \label{eq:89}
  \G_{ij}(n,k,t)=\left\{F\in{[n]\choose k}:|F\cap[t+i+j]|\ge t+i\right\}
\end{equation}
for $0\le i+j\le n-t$; e.g.~$\F_i(n,k,t)=\G_{ii}(n,k,t)$. Observe that:

\begin{prop}
\label{p-5}
$\bigl(\G_{ij}(n,k,t),\G_{ji}(n,l,t)\bigr)$ is cross-$t$-intersecting whenever
$0\le i+j\le n-t$.
\end{prop}

We make the following conjecture, generalizing the
Ahlswede--Khachatrian Theorem (i.e.~that Frankl's General Conjecture
is true, cf.~equation~\eqref{eq:63}).

\begin{conjecture}
\label{j-4}
$\displaystyle{N(n,k,l,t)
=\max_{0\le i+j\le n-t}|\G_{ij}(n,k,t)|\cdot|\G_{ji}(n,l,t)|}$. 
\newline Moreover, up to a permutation of $[n]$, the optimal
cross-$t$-intersecting family is of the form
$\bigl(\G_{ij}(n,k,t),\G_{ji}(n,l,t)\bigr)$ for some $i,j$.
\end{conjecture}

Generalizing $M_0$, we define the maximum size $N_0(n,m_k,m_l,k,l,t)$
of the product of the $m_k$-shade with the $m_l$-shade of a pair of
cross-$t$-intersecting families of $k$-subsets and $l$-subsets of
$[n]$, respectively:
\begin{equation}
  \label{eq:58}
  N_0(n,m_k,m_l,k,l,t)
  =\max_{(\A,\B)\in C(n,k,l,t)}|\shadeto {m_k}(\A)|\cdot|\shadeto {m_l}(\B)|. 
\end{equation}
For purposes of our dichotomy, we were exclusively interested in the numbers
$N_0(2m,m,m,k,k,t)$. Thus we define
\begin{equation}
  \label{eq:45}
  N_1(n,m,k,t)=N_0(n,m,m,k,k,t).
\end{equation}

Corresponding to lemma~\ref{l-7} we have:

\begin{lem}
\label{l-13}
For all $1\le t\le k\le m\le n$,
\begin{enumerate}
\item\label{item:21} $\G_{ij}(n,k,t)=0$ for all $i>k-t$,
\item\label{item:22} $\shadeto m(\G_{ij}(n,k,t))=\G_{ij}(n,m,t)$ 
for all $0\le i\le k-t$ with $i+j\le n-t$.
\end{enumerate}
\end{lem}

Then corresponding to conjecture~\ref{j-1} we have:

\begin{conjecture}
\label{j-5}
$\displaystyle{N_0(n,m_k,m_l,k,l,t)
=\max_{\substack{0\le i\le k-t\\ 0\le j\le l-t\\ i+j\le n-t}}|\G_{ij}(n,m_k,t)|
\cdot|\G_{ji}(n,m_l,t)|}$.
\end{conjecture}

Finally, we arrive at the corresponding asymptotic conjecture.

\begin{conjecture}
\label{co:1}
Assume $k(m)=o(m)$ and $\lim_{m\to\infty}k(m)=\infty$. Suppose that
$\lim_{m\to\infty}\frac{t(m)}{\sqrt{k(m)}}=\infty$.~Then
\begin{equation}
\label{eq:90}
\lim_{m\to\infty}
\frac{\sqrt{N_1(2m,m,k(m),t(m))}}{{2m\choose m}}=0.
\end{equation}
\end{conjecture}

We expect that the argument for corollary~\ref{l-11} will generalize, so that one
can obtain conjecture~\ref{co:1} and a consequence of conjecture~\ref{j-5}.

\bibliographystyle{amsalpha}
\bibliography{database}

\end{document}

%% file: macros.tex
\newcommand{\inv}{^{-1}}

\renewcommand{\div}{\mathbin{/}}

\newcommand{\@shade}[2]{\nabla_{#1#2}}
\newcommand{\@shadow}[2]{\Delta_{#1#2}}
\newcommand{\shadeto}[1]{\@shade{\to}{#1}}
\newcommand{\shadowto}[1]{\@shadow{\to}{#1}}

\newcommand{\@@eq}[2]{=_{#1}^{#2}}
\newcommand{\@@neq}[2]{\ne_{#1}^{#2}}
\newcommand{\@@l}[2]{<_{#1}^{#2}}
\newcommand{\@@nl}[2]{\nless_{#1}^{#2}}
\newcommand{\@@le}[2]{\le_{#1}^{#2}}
\newcommand{\@@nle}[2]{\nleq_{#1}^{#2}}
\newcommand{\@@equiv}[2]{\equiv_{#1}^{#2}}

\newcommand{\@in}[1]{\in_#1}
\newcommand{\@eq}[1]{=_#1}
\newcommand{\@le}[1]{\le_{#1}}
\newcommand{\@l}[1]{<_#1}
\newcommand{\@nle}[1]{\nleq_#1}
\newcommand{\@ge}[1]{\ge_#1}

\newcommand{\u@le}[1]{\le^#1}
\newcommand{\u@l}[1]{<^#1}
\newcommand{\u@eq}[1]{=^#1}

\newcommand{\@lefnt}[1]{\@@le#1*}
\newcommand{\@nlefnt}[1]{\@@nle#1*}
\newcommand{\@eqfnt}[1]{\@@eq#1*}

\newcommand{\len}[1]{\@le{#1}}
\newcommand{\lln}[1]{\@l#1}
\newcommand{\gen}[1]{\@ge{{#1}}}

\newcommand{\str}{\mathrm{str}}

\newcommand{\asym}{\mathrm{asym}}
\newcommand{\trn}{\mathrm{tran}}
\newcommand{\negstr}{-\str}

\newcommand{\strlen}[1]{\@@le{#1}\str}
\newcommand{\nstrlen}[1]{\@@nle{#1}\str}
\newcommand{\strln}[1]{\@@l{#1}\str}
\newcommand{\nstrln}[1]{\@@nl{#1}\str}
\newcommand{\strequiv}[1]{\@@equiv{#1}\str}

\newcommand{\lenn}[2]{\@@le{#1}{#2}}
\newcommand{\nlenn}[2]{\@@nle{#1}{#2}}
\newcommand{\lnn}[2]{\@@l{#1}{#2}}
\newcommand{\nlnn}[2]{\@@nl{#1}{#2}}

\newcommand{\eqnn}[2]{\@@eq{#1}{#2}}
\newcommand{\neqnn}[2]{\@@neq{#1}{#2}}
\newcommand{\asymle}{\u@le\asym}
\newcommand{\asyml}{\u@l\asym}

\newcommand{\strtrnlen}[1]{\@@le{#1}{\str,\trn}}
\newcommand{\strtrnln}[1]{\@@l{#1}{\str,\trn}}
\newcommand{\trnle}{\u@le{\trn}}
\newcommand{\trnl}{\u@l{\trn}}
\newcommand{\trneq}{\u@eq{\trn}}
\newcommand{\negstrlen}[1]{\@@le{#1}\negstr}
\newcommand{\nnegstrlen}[1]{\@@nle{#1}\negstr}
\newcommand{\nnegstrln}[1]{\@@nl{#1}\negstr}
\newcommand{\negstreq}[1]{\@@eq{#1}\negstr}

\renewcommand{\aa}{\mathrm{aa}}
\newcommand{\eqaa}{\@eq\aa}
\newcommand{\inaa}{\@in\aa}

\newcommand{\leaa}{\@le\aa}
\newcommand{\nleaa}{\@nle\aa}
\newcommand{\geaa}{\@ge\aa}
\newcommand{\lefntaa}{\@lefnt\aa}
\newcommand{\nlefntaa}{\@nlefnt\aa}
\newcommand{\eqfntaa}{\@eqfnt\aa}

\newcommand{\naa}{\mathrm{naa}}
\newcommand{\eqnaa}{\@eq\naa}
\newcommand{\eqfntnaa}{\@eqfnt\naa}
\newcommand{\lenaa}{\@le\naa}
\newcommand{\lefntnaa}{\@lefnt\naa}

\newcommand{\eqae}{\@eq\ae}
\newcommand{\leae}{\@le\ae}

\newcommand{\Iff}{\espc\mathrm{iff}\espc}

\renewcommand{\and}{\hespc\mathrm{and}\hespc}

\newcommand{\Or}{\hespc\mathrm{or}\hespc}

\newcommand{\two}{\{0,1\}}

\newcommand{\N}{\mathbb N}

\newcommand{\twoseq}{\two^{<\mathbb N}}

\DeclareFontFamily{U}{cmsy}{}
\DeclareFontShape{U}{cmsy}{m}{n}{<12> sfixed * [10] cmsy10 
<10> <9> <8> <7> <6> <5> sfixed * [10] cmsy10}{}
\DeclareSymbolFont{customtwo}{U}{cmsy}{m}{n} 
\DeclareMathSymbol{\sctn}{\mathord}{customtwo}{"78}

\DeclareFontFamily{U}{cmmi}{}
\DeclareFontShape{U}{cmmi}{m}{n}{<20> sfixed * [11] cmmib10 <12> sfixed * [10]
cmmi10 <10> <9> <8> sfixed * [6] cmmi6 <5> <6> <7> sfixed * [5] cmmi5}{}
\DeclareSymbolFont{custom}{U}{cmmi}{m}{n}
\DeclareMathSymbol{\rharpoon}{\mathord}{custom}{"2A}
\newlength{\widt}
\newlength{\widttwo}
\newlength{\hgt}

\newcommand{\espc}{\quad}
\newlength{\@@hespc}
\newcommand{\hespc}{\hspace{\@@hespc}}

\newcommand{\Th}{{^{\mathrm{th}}}}
\newcommand{\der}[1]{\,\mathrm{d}#1}
\renewcommand{\ae}{\mathrm{ae}}

\newcommand{\tu}{\textup}

\newcommand{\A}{\mathcal A}
\newcommand{\B}{\mathcal B}

\newcommand{\F}{\mathcal F}
\newcommand{\G}{\mathcal G}

\newcommand{\hspb}{\mspace{-1.5mu}}
\newcommand{\spc}{\,\,\,}

\newcounter{saveenumi}




%% file: asymptotic.shade.ub.bbl
\providecommand{\bysame}{\leavevmode\hbox to3em{\hrulefill}\thinspace}
\providecommand{\MR}[1]{}
\providecommand{\MRhref}[2]{%
  \href{http://www.ams.org/mathscinet-getitem?mr=#1}{#2}
}
\begin{thebibliography}{EKR61}

\bibitem[AK97]{MR1429238}
Rudolf Ahlswede and Levon~H. Khachatrian, \emph{The complete intersection
  theorem for systems of finite sets}, European J. Combin. \textbf{18} (1997),
  no.~2, 125--136. 

\bibitem[And02]{MR2003b:05001}
Ian Anderson, \emph{Combinatorics of finite sets}, Dover Publications Inc.,
  Mineola, NY, 2002, Corrected reprint of the 1989 edition. 

\bibitem[EKR61]{MR0140419}
P.~Erd{\H{o}}s, Chao Ko, and R.~Rado, \emph{Intersection theorems for systems
  of finite sets}, Quart. J. Math. Oxford Ser. (2) \textbf{12} (1961),
  313--320. 

\bibitem[Eng97]{MR1429390}
Konrad Engel, \emph{Sperner theory}, Encyclopedia of Mathematics and its
  Applications, vol.~65, Cambridge University Press, Cambridge, 1997.

\bibitem[Fra78]{MR519277}
P.~Frankl, \emph{The {E}rd{\H o}s-{K}o-{R}ado theorem is true for {$n=ckt$}},
  Combinatorics (Proc. Fifth Hungarian Colloq., Keszthely, 1976), Vol. I,
  Colloq. Math. Soc. J\'anos Bolyai, vol.~18, North-Holland, Amsterdam, 1978,
  pp.~365--375. 

\bibitem[Hir08]{Hirs1}
James Hirschorn, \emph{Nonhomogeneous analytic families of trees},
  \href{http://arxiv.org/abs/0807.0147v2}{arXiv:0807.0147v2}, 2008.

\bibitem[Kos89]{MR1044236}
A.~V. Kostochka, \emph{An upper bound on the capacity of the boundary of an
  antichain in an {$n$}-dimensional cube}, Diskret. Mat. \textbf{1} (1989),
  no.~3, 53--61. 

\bibitem[MT89]{MR90g:05008}
Makoto Matsumoto and Norihide Tokushige, \emph{The exact bound in the {E}rd{\H
  o}s-{K}o-{R}ado theorem for cross-intersecting families}, J. Combin. Theory
  Ser. A \textbf{52} (1989), no.~1, 90--97. 

\bibitem[She94]{MR1303493}
Saharon Shelah, \emph{How special are {C}ohen and random forcings, i.e.\
  {B}oolean algebras of the family of subsets of reals modulo meagre or null},
  Israel J. Math. \textbf{88} (1994), no.~1-3, 159--174. 

\bibitem[Usp37]{63.1069.01}
J.~V. Uspensky, \emph{{Introduction to mathematical probability.}}, {IX+ 411 p.
  New York, London, McGraw-Hill Book Co }, 1937 (English).

\end{thebibliography}
